\documentstyle{amsppt}
\documentstyle{psfig}

\mag\magstep1
\nologo

\define\z{\Bbb Z}
\define\q{\Bbb Q}
\def\r{\Bbb R}
\define\p{\Bbb P}
\define\cc{\Bbb C}

\define\zz{\z/2}

\define\xr{X(\r)}

\define\ox{\Cal O_X}
\define\co{\Cal O}

\define\HH{\Cal H}
\define\slr{\operatorname{SL(2,\r)}}
\define\pslr{\operatorname{PSL(2,\r)}}
\define\slz{\operatorname{SL(2,\z)}}
\define\psl{\operatorname{PSL(2,\z)}}
\define\mat{\left(\matrix a&b\\ c&d\endmatrix\right)}
\define\smat{\left(\matrix a&-b\\ -c&d\endmatrix\right)}

\topmatter
\title Surfaces elliptiques r\'eelles et in\'egalit\'e
de Ragsdale-Viro
\endtitle
\author Fr\'ed\'eric Mangolte\endauthor


\affil mangolte\@univ-savoie.fr \\
Fax: (33) 4 79 75 87 42\\  \\
English title: Real elliptic surfaces and
Ragsdale-Viro inequality\endaffil

\address Laboratoire de Math\'ematiques, Universit\'e
de Savoie, F-73~376 Le Bourget du Lac Cedex, France, Tel: (33) 4 79 75 86 60,
Fax: (33) 4 79 75 87 42\endaddress
\email mangolte\@univ-savoie.fr\endemail

\abstract
On a real regular elliptic surface without multiple fiber, the
Betti number $h_1$ and the Hodge number $h^{1,1}$ are related by
$h_1\leq h^{1,1}$. We prove that it's always possible to deform
such algebraic surface to obtain $h_1=h^{1,1}$. Furthermore, we can
impose that each homology class can be represented by a real
algebraic curve. We use a real version of the modular construction of
elliptic surfaces.
\endabstract

\endtopmatter

\hskip 12pt
\vbox{\advance\hsize-2\indenti
  \eightpoint \noindent
  {\smc R\'esum\'e.\enspace}Pour une surface elliptique r\'eelle r\'eguli\`ere et sans fibre multiple, le nombre de
Betti $h_1$ et le nombre de Hodge $h^{1,1}$ sont li\'es par
l'in\'egalit\'e $h_1\leq h^{1,1}$. On montre qu'on peut toujours d\'eformer une telle
surface alg\'ebrique pour obtenir $h_1=h^{1,1}$. De plus, on peut
imposer que chaque classe d'homologie soit repr\'esentable par une
courbe alg\'ebrique r\'eelle. On utilise une version adapt\'ee au cas r\'eel
de la construction des surfaces elliptiques modulaires.
\par\unskip}

\document

\head  1. Introduction \endhead
\medskip
Une vari\'et\'e projective (ou quasi-projective) $X$ est dite r\'eelle ou d\'efinie sur
$\r$ lorsque $X$ est munie d'une involution anti-holomorphe $\sigma_X$ appel\'ee
structure r\'eelle. Dans ce cas, on appelle partie r\'eelle de $X$  et on note $\xr$,
l'ensemble des points fixes de $\sigma_X$. Par convention, une surface est
projective lisse. Une surface est dite elliptique s'il existe une application
holomorphe surjective $\pi\: X\to \Delta$ o\`u $\Delta$ est une courbe lisse
compacte et la fibre g\'en\'erique de $\pi$ est une courbe de genre 1. Une telle
surface est dite elliptique r\'eelle lorsque $X$ et $\Delta$ sont d\'efinies sur
$\r$ et $\pi\circ\sigma_X=\sigma\circ\pi$ o\`u $\sigma$ est la structure r\'eelle de
$\Delta$.

D'apr\`es un r\'esultat de V.~Kharlamov, [Kha], une
surface elliptique r\'eelle $X$ r\'eguli\`ere (i.e. telle que $H^1(X,\co_X)=\{0\}$) et
sans fibre multiple,  v\'erifie l'in\'egalit\'e de Ragsdale-Viro
$$
h_1(\xr)\leq h^{1,1}(X)\tag 1.1
$$
o\`u $h_1(\xr)$ est le rang du premier groupe d'homologie $H_1(\xr,\zz)$ de
$\xr$ et $h^{1,1}(X)=\dim_\cc H^1(X,\Omega_X^1)$. Maintenant, on note
$h^{\roman{alg}}_1(X(\r))$ le rang du sous-groupe de $H_1(\xr,\zz)$
engendr\'e par les classes fondamentales de courbes al\-g\'eb\-ri\-ques
r\'eelles (voir [BH]). On a alors
$h^{\roman{alg}}_1(X(\r))\leq h_1(\xr)$.


Le nombre $h^{1,1}$ est invariant par d\'eformation et le r\'esultat
principal de ce travail est~:

\proclaim{(1.2) Th\'eor\`eme}Toute surface elliptique r\'eguli\`ere sans fibre multiple
peut \^etre d\'eform\'ee sur $\cc$ en une surface elliptique r\'eelle telle que
$h^{\roman{alg}}_1=h_1=h^{1,1}$.\endproclaim

\remark{Remarques}
Une surface d'Enriques $X$ est elliptique r\'eguli\`ere avec deux fibres
multiples. Pour une telle surface, on a $h^{1,1}(X)=10$, mais il  n'existe aucune
surface d'Enriques r\'eelle telle que $h^{\roman{alg}}_1=h_1=10$. En revanche,
il existe des surfaces d'Enriques r\'eelles telles que $h_1^{\roman{alg}}=11$ et
$h_1=12$, cf. [DKh] et [MavH].

On ne peux pas esp\'erer une version du th\'eor\`eme~1.2 avec d\'eformation sur
$\r$. En effet, m\^eme si $h_1(\xr)=h^{1,1}(X)$, il peut n'exister aucune
d\'eformation sur $\r$ de $X$ telle que $h^{\roman{alg}}_1=h^{1,1}$, voir \`a ce sujet
le comportement des surfaces K3, [Ma].
\endremark\medskip

Le point cl\'e dans la preuve du th\'eor\`eme~1.2 est la
construction d'une suite $\{X_k\}_{k\geq 1}$ de surfaces elliptiques r\'eelles qui
v\'erifient $h^{\roman{alg}}_1(X_k(\r))=h^{1,1}(X_k)$, cf. th\'eor\`eme~5.3.

Les surfaces $X_k$ sont des surfaces elliptiques modulaires au sens de T.
Shioda [Sho] dont on a adapt\'e les constructions au cas r\'eel. Jusqu'\`a pr\'esent, les
surfaces elliptiques modulaires ne semblent pas avoir \'et\'e utilis\'ees en g\'eom\'etrie
alg\'ebrique r\'eelle.

\bigskip
Le plan de cet article est organis\'e comme suit~: la section~2 est consacr\'ee aux
pr\'eliminaires sur les courbes de genre 1 r\'eelles, les r\'ef\'erences sont [Si3] et [Sil].
Dans la section~3, on montre comment utiliser les constructions analytiques de
K.~Kodaira [Ko] en r\'eel. Gr\^ace \`a la classification des fibres singuli\`eres des
pinceaux r\'eels de courbes elliptiques due \`a R.~Silhol [Si1], on r\'eduit alors le
probl\`eme \`a la construction de surfaces elliptiques avec fibres singuli\`eres
donn\'ees a priori. La section~4 est consacr\'ee aux surfaces elliptiques modulaires
r\'eelles. On adapte les travaux de T.~Shioda [Sho] et M.~Nori [No] pour r\'eduire le
probl\`eme \`a la recherche de groupes Fuchsiens arithm\'etiques particuliers.
Enfin, en section~5, on donne les domaines fondamentaux d'une suite de
groupes $\Gamma_k$ qui servent de base \`a la construction des surfaces $X_k$.

Je tiens \`a remercier V.~Kharlamov pour m'avoir indiqu\'e le r\'esultat (1.1) et
P.~Schmutz~Schaller pour son aide dans la construction des domaines
fondamentaux de la section~5.
\bigskip
\head 2. Pr\'eliminaires : Courbes de genre 1 r\'eelles \endhead
\medskip
On note $\HH=\{z\in\cc/\ \Im(z)>0\}$ le demi-plan sup\'erieur.
L'action du groupe $\pslr=\slr/\{\pm 1\}$ sur $\HH$ est not\'ee $z\mapsto A.z$ o\`u
$A.z=\frac{az+b}{cz+d}$ si $A$ est repr\'esent\'e par $\mat$, $ad-bc=1$.
L'involution $\sigma_\HH\: z\mapsto -\bar z$ de $\HH$ est anti-holomorphe.
On note $S\: \slr\to\slr$, $\mat\mapsto\smat$, clairement $S$ induit une
involution sur $\pslr$ que l'on notera encore $S$.  Soit $z\in\HH$, et soit
$A\in\pslr$, alors $\sigma_\HH(A.\sigma_\HH(z))=S(A).z$.

Soit $\Gamma$ un groupe Fuchsien (i.e. un sous-groupe discret de $\pslr$),
$\sigma_\HH$ induit une structure r\'eelle sur le quotient $\HH/\Gamma$ si et
seulement si $\sigma_\HH\Gamma=\Gamma\sigma_\HH$ i.e. si et seulement si
$\Gamma$ est stable par $S$.
Dans toute la suite, ce sont les sous-groupes d'indice fini du groupe modulaire
$\psl\subset\pslr$ qui vont nous int\'eresser.

Soit $C$ une courbe projective lisse de genre 1, alors il existe $\tau\in \HH$ tel
que $C=\cc/(\z+\tau\z)$.
R\'eciproquement, soit $\tau\in\HH$, on note $C(\tau)$ la courbe quotient
$\cc/(\z+\tau\z)$.  La fonction modulaire elliptique est not\'ee $j \: \HH\to \cc$. Le
signe $\cong$ signifie isomorphe sur $\cc$. La proposition suivante est classique.

\proclaim{(2.1) Proposition} Soit $(\tau,\tau')\in\HH\times\HH$, les
trois assertions suivantes sont \'equi\-va\-lentes.
\roster
\item $C(\tau)\cong C(\tau')$ ;
\item  $j(\tau)=j(\tau')$ ;
\item il existe $A\in \psl$, tel que $A.\tau=\tau'$.
\endroster
\endproclaim

On note donc $j(C)=j(\tau)$ pour $\tau$ quelconque tel que $C(\tau)\cong C$.

La fonction $j \: \HH\to \cc$ est d\'efinie sur $\r$, c'est-\`a-dire que $\forall
\tau\in\HH$, $\overline{j(\tau)}= j(\sigma_\HH( \tau))$.

\proclaim{(2.2) Proposition} Soit $C$ une courbe projective lisse de genre 1, les
trois assertions suivantes sont \'equivalentes.
\roster
\item $C$ peut \^etre d\'efinie sur $\r$ ;
\item $j(C)\in\r$ ;
\item $\exists\tau\in\HH/\   C(\tau)\cong C$  et $2\Re(\tau)\in\z$.
\endroster
\endproclaim

Soit $\tau\in\HH$ v\'erifiant $2\Re(\tau)\in\z$, la structure r\'eelle induite sur
$C(\tau)=\cc/(\z+\tau\z)$ par la conjugaison complexe de $\cc$ est telle que
$C(\tau)(\r)\ne\emptyset$.

\proclaim{(2.3) Proposition}
Lorsque $C$ est munie d'une structure r\'eelle, on a

si ${}^\#C(\r)= 2$, alors $j(C)\geq 1$ et $\exists \tau\in\HH/\ C(\tau)\cong C$
et $\Re(\tau) \in\z$

si ${}^\#C(\r)=1$, alors $j(C)\leq 1$ et $\exists \tau\in\HH/\  C(\tau)\cong
C$  et $2\Re(\tau)$ est un entier impair.
\endproclaim

\remark{(2.4) Remarque}Le cas $j(C)=1$ correspond aux courbes $y^2=x^3+x$ et
$y^2=x^3-x$ qui sont isomorphes sur $\Bbb C$ par $(x,y)\mapsto (ix,\zeta)$, o\`u
$\zeta^2=-i$. La premi\`ere courbe poss\`ede une partie r\'eelle connexe, la
partie r\'eelle de la seconde courbe poss\`ede deux composantes connexes.
\endremark
\bigskip
\head 3. Surfaces elliptiques r\'eelles \endhead
\medskip
On consid\`ere une  surface elliptique r\'eelle $\pi\: X \to \Delta$. On suppose
que $\pi$ poss\`ede au moins une fibre singuli\`ere, est sans fibre multiple et
qu'aucune courbe exceptionnelle n'est contenue dans une fibre.

Par hypoth\`ese, $\Delta$ est munie d'une involution anti-holomorphe $\sigma$.
On note $C_u=\pi^{-1}\{u\}$ la fibre au dessus de $u\in\Delta$. On consid\`ere un
ensemble fini $\Sigma\subset\Delta$ stable par
$\sigma$ tel que $C_u$ est lisse pour tout $u\in\Delta'=\Delta\setminus\Sigma$.
Comme $\Sigma$ est stable par $\sigma$, $\Delta'$ munie de la restriction de
$\sigma$ est encore r\'eelle.

L'invariant fonctionnel de $\pi$ est la fonction m\'eromorphe $J\: \Delta'\to\cc$,
$u\mapsto j(C_u)$ (cf. Section~2). Par construction $J$ est d\'efinie sur $\r$, i.e.
$J\circ\sigma=\bar J$. On peut prolonger $J$ en une fonction
holomorphe $J\:\Delta\to\p^1(\cc)$ qui v\'erifie $J\circ\sigma=conj\circ J$ o\`u
$conj$ est la conjugaison complexe sur $\p^1(\cc)$.

Maintenant $X_{\vert_{\Delta'}}=\pi^{-1}(\Delta')\to\Delta'$ est un
fibr\'e diff\'erentiel en tores, donc les groupes d'homologie des fibres
$\{H_1(C_u(\cc),\z)\}_{u\in\Delta'}$ forment un faisceau localement constant au
dessus de $\Delta'$. On peut \'etendre ce faisceau \`a $\Delta$, [Ko, \S 7]. Ce
faisceau \'etendu $G$ est l'invariant homologique de $\pi$. De m\^eme que
$J$, $G$ est d\'efini sur $\r$.

Si $J$ est non constante, on peut \'etendre l'ensemble fini $\Sigma$, avec
$\sigma(\Sigma)=\Sigma$, pour obtenir $\forall u\in\Delta'$,
$J(u)\not\in\{0,1,\infty\}$. On note
$\rho'\:\pi_1(\Delta')\to\slz$ l'homo\-mor\-phis\-me de monodromie associ\'e. Le
morphisme
$\rho'$ est une repr\'esentation de
$\pi_1(\Delta')$ qui d\'etermine et est d\'etermin\'e par le faisceau $G$.
Comme
$\pi$ est r\'eelle, pour tout $\alpha\in
\pi_1(\Delta')$, on a $\rho'(\sigma_*(\alpha))=S(\rho'(\alpha))$ o\`u
$\sigma^*$ est l'involution induite par $\sigma$ sur $\pi_1(\Delta')$ et $S$ est
d\'efinie au d\'ebut de la section 2.

Soit $u\in\Sigma$ et soit $\alpha\in\pi_1(\Delta')$ l'\'el\'ement repr\'esent\'e par un
lacet simple tournant dans le sens positif autour de $u$.
Le point $u$ est un p\^ole de $J$ si et seulement si l'ordre de $\rho'(\alpha)$ est
infini. Dans ce cas, $\rho'(\alpha)$ est conjugu\'ee dans $\slz$ \`a une matrice d'une
des formes
$$
\left(\matrix 1&m\\ 0&1\endmatrix\right)\text{ ou }\left(\matrix
-1&-m\\ 0&-1\endmatrix\right)
\text{ avec } m>0\tag 3.1
$$
cf. [Ko \S 9]. Dans le premier cas, on
dit que la fibre $C_{u}$ est du type $I_{m}$ et dans le second cas, du type $I^*_{m}$.

Si $u$ est un p\^ole de $J$ appartenant \`a $\Delta(\r)$, $C_u$ est r\'eelle et les types
r\'eels possibles de $C_u(\r)$ sont classifi\'es par la table (3.2) ci-dessous extraite de
[Si2, Th. VII(1.5)]. On suppose que $\pi$ admet une section r\'eelle, soit $u'$ un
point de
$\Delta'(\r)$ voisin de $u$. Alors la fibre $C_{u'}$ est lisse et $C_{u'}(\r)$ poss\`ede
une ou deux composantes connexes. Dans la deuxi\`eme colonne de la table, on a
indiqu\'e le nombre de composantes connexes de $C_{u'}(\r)$ s'il reste constant
lorsque $u'$ varie au voisinage de $u$ dans $\Delta'(\r)$. On a indiqu\'e "$*$" si
${}^\#C_{u'}(\r)$ change au voisinage de $u$. Les colonnes 3 et 4 donnent la
charact\'eristique d'Euler topologique de $C_u(\r)$ et $C_u$. La derni\`ere colonne
donne le nombre
$\xi(C_u)-1$ de composantes irr\'eductibles de $C_u$ qui ne rencontrent pas la
section.

$$
\boxed{
\left.\matrix
\\ \text{Type complexe}\\  \\ \matrix\format \l\\
 I_m^*, m\text{ pair}\\  \\   I_m^*, m\text{ impair}\\  \\  I_m, m\text{ pair}\\ \\
\\  \\    I_m, m\text{ impair}\\ \\
\endmatrix
\endmatrix
\right\vert
\left.\matrix
\\ {}^\#C_{u'}(\r)\\ \\ \matrix
2 \\ 1\\  *\\ * \\   2\\  2\\  1\\  1\\  * \\  *\\
\endmatrix\endmatrix
\right\vert
\left.\matrix
\\  \chi(C_{u}(\Bbb R))\\ \\ \matrix
-m-4  \\ -m-2\\ -m-4\\  -m-2\\ -m \\ 0 \\ -m\\ 0\\ -m \\  1
\endmatrix\endmatrix
\right\vert
\left.\matrix
\\  \chi(C_{u})\\ \\ \matrix
m+6 \\  \\  \\   \\  m \\  \\  \\   \\   \\ \\
\endmatrix\endmatrix
\right\vert
\left.\matrix
\\  \xi(C_u)-1\\ \\ \matrix
m+4 \\  \\  \\   \\  m-1 \\  \\  \\   \\   \\ \\
\endmatrix\endmatrix\right.
}\tag 3.2
$$

Pour simplifier, et comme cela suffira pour la suite, on a laiss\'e de c\^ot\'e les fibres
singuli\`eres qui pourraient appara\^{\i}tre en dehors des p\^oles de $J$.
Remarquons que $\deg(J)=\sum_{l=1}^t m_l$ o\`u $t$ est le nombre de p\^oles de $J$.
Pour une surface elliptique, la charact\'eristique d'Euler topologique vaut 12 fois
la charact\'eristique d'Euler holomorphe $\chi(\ox)$. Si toutes les fibres singuli\`eres
de
$\pi$ sont de type
$I$ ou
$I^*$, on tire de l'avant-derni\`ere colonne de la table (3.2)~:
$$
12\chi(\ox)=\mu+6\nu(I^*)\tag 3.3
$$
o\`u $\mu= \deg(J)$ et $\nu(I^*)$ est le nombre de fibres du type $I^*$.

Au vu de la table~3.2, on n'a pas de formule aussi g\'en\'erale pour les invariants
r\'eels, voir l'exemple~4.6 et la preuve du th\'eor\`eme~5.3.
\bigskip
\head 4. Surfaces elliptiques modulaires r\'eelles \endhead
\medskip
Soit $\Gamma\subset\psl$ un sous-groupe d'indice fini, en particulier $\Gamma$
est un groupe Fuchsien arithm\'etique. Comme $\Gamma$ op\`ere sur $\p^1(\q)$,
le quotient
$$
\Delta_\Gamma=(\HH\cup\p^1(\q))/\Gamma
$$
est bien d\'efini. Par hypoth\`ese, le nombre de classes de conjugaisons
paraboliques de $\Gamma$ est fini donc $\Delta_\Gamma$ est une courbe
projective lisse cf. [Shm, \S 1.3 et \S 1.5]. Si $\Gamma\subset\Gamma'$ est un
sous-groupe, alors l'application canonique de $\HH/\Gamma$ sur
$\HH/\Gamma'$ s'\'etend en une application holomorphe de $\Delta_\Gamma$ sur
$\Delta_{\Gamma'}$. En particulier, on a une application holomorphe
$$
J_\Gamma\:\Delta_\Gamma\to\p^1
$$
en prenant $\Gamma'=\psl$ et en identifiant $\Delta_{\Gamma'}$ avec $\p^1$ via
la fonction $j$, cf. section~2.
Si on suppose de plus que $\Gamma$ est stable par $S$, $\Delta_\Gamma$ et
$J_\Gamma$ sont naturellement d\'efinis sur $\r$.

L'application $J_\Gamma$ est ramifi\'ee seulement au dessus des trois
points 0, 1 et $\infty$, cf. [No, prop. 2.1]. De plus, le degr\'e de $J_\Gamma$ est
\'egal \`a l'indice de $\Gamma$ dans $\psl$
$$
\deg J_\Gamma=[\psl : \Gamma]
$$

Notons $\Delta'=\Delta_\Gamma\setminus J^{-1}_\Gamma\{0,1,\infty\}$ et
$p\:\HH\to\Delta'$ le rev\^etement universel. Par construction, il existe un
isomorphisme $w\:\HH\to\HH$ qui fait commuter le diagramme
$$
\CD
\HH @>w>>\HH\\
@VpVV  @VjVV\\
\Delta' @>J_\Gamma>> \p^1{\ssize\setminus\{0,1,\infty\}}
\endCD\tag 4.1
$$
et qui v\'erifie $w\circ\sigma_\HH=\sigma_\HH\circ w$.
De l\`a, il existe une unique repr\'esentation
$$
\rho\:\pi_1(\Delta')\to\Gamma\subset\psl\;.
$$
telle que $\forall\alpha\in\pi_1(\Delta')$, $\forall \tilde u\in \HH$,
$\rho(\alpha).w(\tilde u)=w(\alpha.\tilde u)$
et $\rho(\sigma_*(\alpha))=S(\rho(\alpha))$.
Maintenant chaque rel\`evement
$\rho'\:\pi_1(\Delta')\to\slz$ de $\rho$ v\'erifie encore
$\rho'(\sigma_*(\alpha))=S(\rho(\alpha))$.
Il est facile d'adapter au cas r\'eel la construction de Kodaira  [Ko, p. 578],
pour montrer qu'il existe une surface elliptique
$$
\pi\: X\to\Delta_\Gamma
$$
d\'efinie sur $\r$ et poss\'edant une section r\'eelle ayant pour invariants
$J_\Gamma$ et $\rho'$.

\definition{(4.2) D\'efinition} Si $-1\not\in\rho'(\pi_1(\Delta'))$, une telle surface $
\pi\: X\to\Delta_\Gamma
$ est
appel\'ee surface elliptique modulaire r\'eelle.\enddefinition

Pour simplifier, et ceci correspond \`a la restriction sur les types de fibres
singuli\`eres d\'ej\`a faite en section 3, on se restreint au cas o\`u
$\Delta_\Gamma$ est sans point elliptique i.e. $\Gamma$ est sans torsion, on a
alors (cf. [Shm] ou [Sho (4.6)])
$$
2g-2+t=\frac 16\mu \tag 4.3
$$
o\`u $g$ est le genre de $\Delta_\Gamma$, $t$ est le nombre de cusps et
$\mu=[\psl :\Gamma]$.

On note $\{u_l\}=J^{-1}_\Gamma\{0,1,\infty\}$ avec pour $1\leq l\leq t$,
$J(u_l)=\infty$. Alors les fibres singuli\`eres de $\pi\:
X\to\Delta_\Gamma$ sont les $t$ fibres $C_{u_l}=\pi^{-1}\{u_l\}$, $l\leq t$ (cf.
[Sho]) et le type de la fibre
$C_{u_l}$ est determin\'e par $\rho'(\alpha_l)$ o\`u $\alpha_l\in\pi_1(\Delta')$ est
l'\'el\'ement repr\'esent\'e par un lacet simple tournant dans le sens positif autour de
$u_l$. Soit $l\leq t$, choisissons un point $z\in\q\cup\{\infty\}$ repr\'esentant $u_l$.
Comme $u_l$  est un cusp, le g\'en\'erateur du stabilisateur de $z$
dans $\Gamma$ est conjugu\'e dans
$\psl$ \`a
$\left(\matrix 1&m_l\\ 0&1\endmatrix\right)\mod \{\pm 1\}$, $m_l>0$
et la fibre $C_{u_l}$ est du type $I_{m_l}$ ou $I_{m_l}^*$. Les invariants
num\'eriques de $X$ sont alors donn\'es par
$$
\aligned
q(X)&=g(\Delta_\Gamma)\\
12\chi(\ox)&=\mu+6\nu(I^*)
\endaligned \tag 4.4
$$

Sachant que $\mu=\sum_{l=1}^t m_l$ et $\nu(I^*)$ est le
nombre de fibre de type $I^*$.

\proclaim{(4.5) Lemme} Si $\Gamma$ est sans
torsion et $g(\Delta_\Gamma)=0$, il existe exactement $t-1$ rel\`evements
distincts $\rho'\:\pi_1(\Delta')\to\slz$ de $\rho\:\pi_1(\Delta')\to\psl$ qui
v\'erifient
$$
\rho'(\alpha_l)=1\qquad\text{si}\qquad\rho(\alpha_l)=1
$$
o\`u $t$ est le nombre de cusps de $\Delta_\Gamma$.
\endproclaim

C'est un cas particulier de [No, prop. 2.3].

\example{(4.6) Exemple} On consid\`ere le groupe de congruence
$$
\Gamma(2)=\left\{\mat\in\slz/\ a\equiv d\equiv 1 \mod 2, b\equiv c\equiv 0
\mod 2\}\right\}/\{\pm 1\}
$$
d'indice $\mu=6$ dans $\psl$  (cf. e.g. [Kat, chap. V]).

\bigskip
\vbox{\centerline{\psfig{file=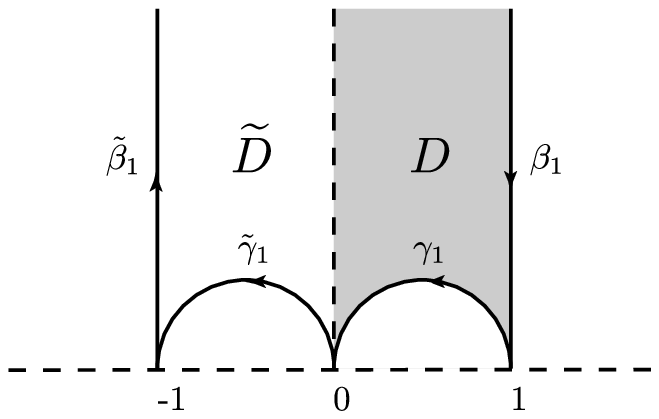}}\centerline{\bf Fig. 1}}
\bigskip

On a repr\'esent\'e sur la figure 1 un domaine fondamental pour $\Gamma(2)$. Il
est limit\'e  par les g\'eod\'esiques $\beta_1,\gamma_1,\tilde
\gamma_1,\tilde\beta_1$, avec les identifications $l\sim \tilde l$ pour une
g\'eod\'esique $l$.  Ce domaine est stable par $\sigma_\HH$. De l\`a, la courbe
quotient $\Delta=(\HH\cup\p^1(\q))/\Gamma(2)$ est munie de la structure
r\'eelle quotient et sa partie r\'eelle $\Delta(\r)$ est repr\'esent\'ee dans $\HH$ par
la r\'eunion des g\'eod\'esiques $\beta_1\cup\gamma_1$ et du demi-axe
imaginaire. En particulier, les trois cusps $u_1,u_2, u_3$ de $\Delta$, repr\'esent\'es
par $0,1$ et $\infty$ appartiennent \`a $\Delta(\r)$.

Remarquons au vu des identifications que $g(\Delta)=0$ i.e.
$\Delta\cong \p^1(\cc)$.

Consid\'erons une surface modulaire $\pi\: X\to \Delta$ associ\'ee \`a un
rel\`evement, cf.(4.1),
$$
\rho'\:\pi_1(\Delta\setminus \{u_l\})\to \slz
$$

La fibration $\pi$ poss\`ede exactement trois fibres singuli\`eres du type $I_m$ ou
$I_m^*$ et d'apr\`es (4.4), on a~:
$$
\aligned
q(X)&=0\\
12\chi(\ox)&=6+6\nu(I^*)
\endaligned\tag 4.7
$$

Ce qui impose que le nombre de fibres de type $I^*$ est impair donc \'egal \`a
1 ou 3. Dans le premier cas, les trois fibres sont $I_2,I_2,I_2^*$. Alors $\chi(\ox)=1$
et $X$ est une surface rationnelle r\'eelle telle que
$h_1^{\roman{alg}}(\xr)=h_1(\xr)=h^{1,1}(X)=10$. Dans le deuxi\`eme cas, les fibres
singuli\`eres sont trois $I_2^*$, $\chi(\ox)=2$ et  $X$ est une surface K3 r\'eelle telle
que $h_1^{\roman{alg}}(\xr)=h_1(\xr)=h^{1,1}(X)=20$.
\endexample

\bigskip
\head 5. Construction de surfaces extr\'emales\endhead
\medskip
Pour une  surface elliptique $X$ r\'eguli\`ere, on a
$$
h^{1,1}(X)=10\chi(\ox)\;.\tag 5.1
$$

En effet, d'apr\`es le th\'eor\`eme de d\'ecomposition de Hodge,
$$
h^{1,1}(X)=\dim_\cc H^2(X,\cc)-2\dim_\cc H^2(X,\ox)\;.
$$

L'hypoth\`ese de r\'egularit\'e
nous donne $\dim_\cc H^2(X,\cc)=12\chi(\ox)-2$ (voir aussi (3.3))  et
$\dim_\cc H^2(X,\ox)=\chi(\ox)-1$ d'o\`u le r\'esultat.

\definition{(5.2) D\'efinition} Une surface elliptique r\'eguli\`ere $\pi\: X\to\p^1$ est
dite normalis\'ee si et seulement si
\roster
\item aucune courbe exceptionnelle n'est contenue dans une fibre ;
\item$\pi$ est sans fibre multiple.
\endroster
\enddefinition

\proclaim{(5.3) Th\'eor\`eme} Pour tout entier $k>0$, il existe une surface
elliptique r\'eelle normalis\'ee  $\pi\: X_k\to \p^1$ qui v\'erifie
$$
\chi(\co_{X_k})=k\quad\text{et}\quad h_1^{\roman{alg}}(X_k(\r))=10k
$$
\endproclaim

Avant de prouver le th\'eor\`eme~5.3, on montre comment on peut d\'eduire de
ce r\'esultat le th\'eor\`eme~1.2 de l'introduction.

\proclaim{(5.4) Th\'eor\`eme (Kodaira)} Toute surface elliptique r\'eguli\`ere $X$
provient d'une surface elliptique normalis\'ee $Y$ par une succession de
$\eta$ \'eclatements de points et de $\delta$ transformations logarithmiques
d'ordres respectifs $n_1\dots,n_\delta$.

Deux surfaces elliptiques r\'eguli\`eres
$X$ et $X'$ sont d\'efor\-ma\-tions l'une de l'autre si et seulement si
\roster
\item $Y$ et $Y'$ sont d\'eformations l'une de l'autre
\item $\eta=\eta'$
\item $\delta=\delta'$ et $n_l=n'_l$, $\forall l\in\{1,\dots,\delta\}$
\endroster
\endproclaim
Cet \'enonc\'e est tir\'e de [Pe, p.306]. Le th\'eor\`eme suivant se trouve en premi\`ere
page de [Kas].

\proclaim{(5.5) Th\'eor\`eme (Kas)} Deux surfaces elliptiques normalis\'ees $Y$ et $Y'$
sont d\'efor\-ma\-tions l'une de l'autre si et seulement si
$$
\chi(\co_{Y})=\chi(\co_{Y'})
$$
\endproclaim

\demo{Preuve du th\'eor\`eme~1.2}
Soit $V$ une surface r\'eelle et $W\to V$ un \'eclatement centr\'e en un point de
$V(\r)$, alors $W$ est une surface r\'eelle et
$$
\align
h^{1,1}(W)&=h^{1,1}(V)+1,\\
h_1^{\roman{alg}}(W(\r))&=h_1^{\roman{alg}}(V(\r))+1\;.
\endalign
$$

Soit $X$ une surface elliptique r\'eguli\`ere, d'apr\`es (5.4), $X$ provient d'une surface
normalis\'ee  $Y$ par une suite de $\eta$ d'\'eclatements et $\delta$
transformations logarithmiques. Si de
plus $X$ est sans fibre multiple, on a $\delta=0$. Posons
$k=\chi(\co_Y)$, et consid\'erons une surface $X_k$ v\'erifiant les hypoth\`eses du
th\'eor\`eme~5.3, alors
$h_1^{\roman{alg}}(X_k(\r))=h^{1,1}(Y)$ et d'apr\`es le th\'eor\`eme de Kas (5.5), $X_k$
est d\'eformation de $Y$. Notons $X'$ la surface obtenue \`a partir de $X_k$ apr\`es
\'eclatement de $\eta$ points de
$X_k(\r)$. Alors $X'$ est une d\'eformation de $X$ par (5.4)
et $h_1^{\roman{alg}}(X'(\r))=h^{1,1}(X')$.
D'apr\`es (1.1), comme $X'$ est une surface elliptique r\'eelle r\'eguli\`ere et sans fibre
multiple, on a
$$
h_1^{\roman{alg}}(X'(\r))\leq h_1(X'(\r))\leq h^{1,1}(X')\;.
$$
d'o\`u le r\'esultat.

\enddemo

\demo{Preuve du  th\'eor\`eme~5.3}

Le cas $k=1$ est trait\'e dans l'exemple~4.6. Soit $k\geq 2$,
\`a partir du domaine fondamental $D\cup\widetilde D$ de $\Gamma(2)$, cf.
Fig.~1, on construit un nouveau domaine obtenu en adjoignant au
triangle $D\cup\widetilde D$ $k-2$ translat\'es de $D$ \`a droite et leurs
sym\'etriques par rapport au demi-axe imaginaire. Le nouveau domaine est
bord\'e par les g\'eod\'esiques $\beta_{k-1},\gamma_{k-1},\gamma_{k-2},\dots
\gamma_1$, $\tilde\gamma_1,\dots,\tilde\gamma_{k-1},\tilde\beta_{k-1}$
identifi\'ees deux \`a deux par $l\sim\tilde l$.

\bigskip
\vbox{\centerline{\psfig{file=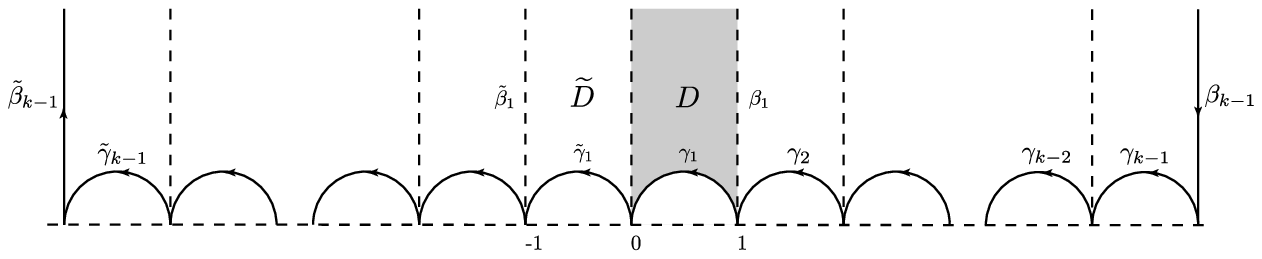}}\centerline{\bf Fig.
2}}
\bigskip

C'est le domaine fondamental d'un sous-groupe
$\Gamma_k\subset\Gamma(2)$ d'indice fini dans $\psl$, cf. e.g. [DR]. Par
construction, la courbe quotient $\Delta_k=(\HH\cup\p^1(\q))/\Gamma_k$ est
sans point elliptique et poss\`ede $k+1$ cusps.

Le domaine fondamental est stable par $\sigma_\HH$ et la courbe
$\Delta_k$ est munie de la structure r\'eelle quotient. Sa partie r\'eelle
$\Delta_k(\r)$ est repr\'esent\'ee dans $\HH$ par la r\'eunion des g\'eod\'esiques
$\beta_{k-1}\cup\gamma_{k-1}\cup\gamma_{k-2}\cup\dots
\gamma_1$ et du demi-axe imaginaire.

Au vu des identifications, on a $g(\Delta_k)=0$ et d'apr\`es (4.3)
$$
\mu(\Gamma_k)=6(k-1)\;.
$$

Suivant la section 4, on a une application $J_{\Gamma_k}\:
\Delta_k\to\p^1$ et on note $\Delta_k'=\Delta_k\setminus
J_{\Gamma_k}^{-1}\{0,1,\infty\}$. On peut construire pour chaque rel\`evement
$$
\rho'\:\pi_1(\Delta_k')\to\slz
$$
de $\rho\:\pi_1(\Delta_k')\to\Gamma_k\subset\psl$ une surface elliptique
modulaire r\'eelle avec section r\'eelle
$$
\pi \:X\to\p^1
$$
d'invariants $(J_{\Gamma_k},\rho')$.
Par construction, $X$ poss\`ede exactement $k+1$ fibres singuli\`eres et ces
fibres sont de type $I_m$ ou $I_m^*$ .

D'apr\`es les formules (4.4), la surface $X$ v\'erifie :
$$
\aligned
q(X)&=0\\
12\chi(\ox)&=6k-6+6\nu(I^*)
\endaligned
$$
o\`u $\nu(I^*)$ est le nombre de fibres de type $I_m^*$.
En particulier, quelque soit le choix de $\rho'$, $X$ est une surface elliptique
r\'eguli\`ere et $h^{1,1}(X)=10\chi(\co_{X})$.

Si $k$ est impair, le nombre de fibres singuli\`eres est pair et d'apr\`es le
lemme~4.5, il existe un rel\`evement $\rho'$ tel que toutes les fibres singuli\`eres
soient de type $I^*$.

Si $k$ est pair, $6k-6$ n'est pas divisible par 12 et
n\'ecessairement l'une des fibres singuli\`ere est de type $I^*$. Donc il existe
encore un rel\`evement de $\rho$ tel que toutes les fibres singuli\`eres soient de
type $I^*$. On peut donc obtenir pour tout $k\geq 2$ une surface  not\'ee $X_k$
dont la liste des fibres singuli\`eres est $I^*_{m_1},\dots,I^*_{m_{k+1}}$ avec $\sum_l
m_l=6k-6$. On a alors
$$
\chi(\co_{X_k})=k, \qquad h^{1,1}(X_k)=10k
$$

Par construction, les $k+1$ cusps appartiennent \`a $\Delta_k(\r)$ donc les
$k+1$ fibres de type $I^*_{m_l}$ sont r\'eelles. Comme
$\Gamma_k\subset\Gamma(2)$, les cusps sont chacun conjugu\'e (dans
$\psl$) \`a un cusp $\left(\matrix 1&m_l\\ 0&1\endmatrix\right)$ avec
$m_l$ pair.

On d\'eduit alors de la table~3.2 que le nombre de composantes connexe de la
partie r\'eelle d'une fibre r\'eelle de $\pi$ est constant au voisinage de chaque fibre
singuli\`ere. Donc toutes les fibres lisses de $X_k(\r)$ ont le m\^eme nombre de
composantes connexe. Par ailleurs, si on note $\Omega$ l'image du demi-axe
imaginaire par l'application canonique
$\HH\to\Delta_k=(\HH\cup\p^1(\q))/\Gamma_k$, on a
$\Omega\subset\Delta_k(\r)$. D'apr\`es les propositions~2.2 et
2.3, on a $J(u)>1$ pour presque tout les points $u$ de $\Omega$ donc
${}^\#C_{u}(\r)=2$. D'apr\`es ce qui pr\'ec\`ede, on a plus g\'en\'eralement
${}^\#C_{u}(\r)=2$ pour toute les fibres r\'eelles lisses de $\pi$.

Par construction, $X_k(\r)$ est connexe et
$h_1(X_k(\r))=2-\sum_l\chi(C_{u_l}(\r))$. De nouveau \`a partir de la table (3.2),
maintenant qu'on connait le nombre de composantes au voisinage de chaque
cusp $u_l$, on d\'eduit $\chi(C_{u_l}(\r))=-m_l-4$. De plus, d'apr\`es [Si2,~VII.1], les
composantes irr\'eductibles r\'eelles d'une fibre de type $I^*$ qui ne rencontrent
pas la section sont ind\'ependantes. De l\`a
$$
-\chi(C_{u_l}(\r))=\xi(C_u)-1
$$
et
$$
h_1^{\roman{alg}}(X_k(\r))=2+\sum_l m_l+4(k+1)=10k
$$
\enddemo

\medskip
La proposition suivante sert \`a d\'eterminer le type topologique
des surfaces $X_k(\r)$.

\proclaim{(5.6) Proposition} Soit $X\to\p^1$ une surface elliptique
r\'eelle normalis\'ee,

si $\chi(\ox)$ est pair et $\xr\ne\emptyset$, alors
$$
\xr \text{ est orientable.}
$$

si $\chi(\ox)$ est impair et $X$ admet une section d\'efinie sur $\r$, alors
$$
\xr\ne\emptyset \text{ et $\xr$ est non orientable.}
$$
\endproclaim
\demo{Preuve} Sous les hypoth\`eses consid\'er\'ees, un diviseur
canonique de $X$ est donn\'e par (cf. e.g. [BPV, p. 162])~:
$$
K_X=(\chi(\ox)-2)F
$$
o\`u $F$ est une fibre quelconque. La classe  $w_1(\xr)\in H_1(\xr,\zz)$, duale de
Poincar\'e de la 1\`ere classe de Stiefel-Whitney, est repr\'esent\'ee par $K_X(\r)$ et
pour tout diviseur r\'eel $D$ de $X$, on a
$$
D.K_X\equiv D(\r)\cap K_X(\r)\mod 2\tag 5.7
$$

Dans le cas o\`u $\chi(\ox)$ est pair, $K_X\equiv 0\mod 2$ d'o\`u
$w_1(\xr)=0$ et $\xr$ est orientable.
Maintenant, supposons que $\chi(\ox)$ est impair et notons $S$ la
courbe image dans
$X$ d'une section d\'efinie sur $\r$, on a
$$
S.K_X\equiv 1 \mod 2
$$
donc $\xr\ne\emptyset$ car n\'ecessairement il existe un point de $S\cap
K_X$ dans $\xr$. On a alors $w_1(\xr)\ne 0$ d'apr\`es (5.7) et $\xr$
est non orientable.
\enddemo

\medskip
On note $S_g$ la surface (topologique) orientable de genre $g$ et
$V_q=\#_q\p^2(\r)$ la surface non orientable de charact\'eristique d'Euler $2-q$.

\proclaim{(5.8) Corollaire} Soit $X_k$ une surface v\'erifiant les
hypoth\`eses du th\'eor\`eme~5.3, alors si
$k$ est pair, $X_k(\r)=S_{5k}$ et si $k$ est impair, $X_k(\r)=V_{10k}$.
\endproclaim

Par hypoth\`ese, on a $h_1(X_k(\r))=10k$ et l'orientabilit\'e est donn\'ee
par la proposition~5.6. Il reste \`a montrer que la surface $X_k(\r)$
est connexe. C'est une cons\'equence de l'in\'egalit\'e de Comessatti
$$
2{}^\#\xr-h_1(\xr)\leq h^{1,1}(X)-2(r(X)-1)
$$
valable pour une surface r\'eelle $X$ quelconque o\`u $r(X)$ d\'esigne le
nombre de classes r\'eelles ind\'ependantes dans le groupe de
N\'eron-Severi de $X$. Comme $X_k$ est normalis\'ee, on a
$h_1^{\roman{alg}}(X_k(\r))\leq r(X_k)\leq h^{1,1}(X_k)$, voir par
exemple [Si2, III(1.10)], donc ${}^\#X_k(\r)=1$.


\Refs\nofrills{R\'eferences Bibliographiques}
\widestnumber\key{MavH}

\ref\key BPV \by W. Barth, C. Peters, A. Van de Ven \book Compact complex
surfaces \bookinfo Ergebnisse der Mathematik \publ Springer
\publaddr Berlin Heidelberg \yr 1984 \endref

\ref\key BH \by E. Borel, A. Haefliger \paper La classe d'homologie
fondamentale d'un espace analytique \jour Bul. Soc. Math. France \vol
83 \pages 461--513 \yr 1961\endref

\ref\key DKh \by  A.~Degtyarev, V.~Kharlamov \paper Topological
classification of real Enriques surfaces \jour Topology\vol 35\pages
711--730\issue 3
\yr 1996 \endref

\ref\key DR \by G. de Rham\paper Sur les polygones g\'en\'erateurs de
groupes Fuchsiens\jour L'ens. Math. \vol 17 \issue 1 \yr 1971\pages
49--61\endref

\ref\key Kas \by A. Kas \paper On the deformation types of regular
elliptic surfaces \inbook Complex analysis and algebraic geometry \eds
W. Baily, T. Shioda \publ Cambridge university press \publaddr
Cambridge London New York\yr 1977 \endref

\ref\key Kat \by S. Katok  \book Fuchsian groups \bookinfo Chicago
Lectures in Mathematics \publaddr Chicago
\publ The University of Chicago Press \yr 1992\endref

\ref\key Kha \by V. Kharlamov\paper Communication priv\'ee\year
1997 \endref

\ref\key Ko \by K. Kodaira \paper On the structure of compact
complex analytic surfaces II \& III\jour Amer. J. Math. \vol 88 \& 90\yr 1966
\&1968 \pages 682--721 \& 56--83 \endref

\ref\key Ma \by F. Mangolte \paper Cycles alg\'ebriques sur les
surface K3 r\'eelles \jour Math. Z. \yr 1997 \vol 225 \issue 4 \pages
559--576\endref

\ref\key MavH\by F. Mangolte, J. van Hamel \paper Algebraic cycles and
topology of real Enriques surfaces \jour Compos. Math.\vol 110 \issue 2\pages
215--237\yr 1998\endref

\ref\key No \by M. Nori \paper On certain elliptic surfaces with
maximal Picard number \jour Topology \vol 24 \yr 1984 \pages 175--186
\endref

\ref\key Pe \by U. Persson \paper Horikawa surfaces with maximal
Picard number \jour Math. Ann.\vol 259 \yr 1982 \pages 287--312
\endref

\ref\key Shm \by G. Shimura \book Introduction to the arithmetic theory of
automorphic functions \publ Princeton Univ. Press \yr 1971\endref

\ref\key Sho \by T. Shioda \paper On elliptic modular surfaces
\jour J. Math. Soc. Japan \vol 24 \yr 1972 \pages 20--59 \endref

\ref\key Si1 \by R. Silhol \paper Real algebraic surfaces with
rational or elliptic fiberings \jour Math. Z. \vol 186 \pages 465--499
\yr 1984 \endref

\ref\key Si2 \by R. Silhol \book Real algebraic surfaces
\bookinfo  Lectures notes in Math. \publ Springer \vol 1392 \publaddr
Berlin Heidelberg New York  \yr 1989 \endref

\ref\key Si3 \by R. Silhol \paper Compactifications of moduli spaces
in real algebraic geometry  \jour Invent. Math. \vol 107
\pages151--202 \yr 1992 \endref

\ref\key Sil \by J. H. Silverman \book The arithmetic of elliptic curves
\publ Springer-Verlag \bookinfo Graduate Texts in Math.\vol 106\yr
1986\publaddr New York\endref

\endRefs
\enddocument